\documentclass{article}

\usepackage{PRIMEarxiv}

\usepackage[utf8]{inputenc} 
\usepackage[T1]{fontenc}    
\usepackage{hyperref}       
\usepackage{url}            
\usepackage{booktabs}       
\usepackage{amsfonts}       
\usepackage{nicefrac}       
\usepackage{microtype}      
\usepackage{xcolor}         
\usepackage{graphicx}
\usepackage{subfigure}
\usepackage{amsmath}
\usepackage{amssymb}
\usepackage{mathtools}
\usepackage{amsthm}
\usepackage{placeins}
\usepackage{multirow}
\usepackage[capitalize,noabbrev]{cleveref}
\usepackage{dsfont}
\usepackage{tikz}
\usetikzlibrary{arrows, arrows.meta,chains,positioning,shapes.geometric}
\usepackage{enumerate}
\usepackage{empheq}
\usepackage{tabularx}
\usepackage{comment}
\usepackage{pifont}
\usepackage{wrapfig}
\usepackage{soul} 
\usepackage{caption}
\usepackage{stmaryrd}
\usepackage{bm}
\usepackage{booktabs}
\usepackage{import}

\begin{document}

\title{Learning large scale industrial physics simulations}

\author{Fabien Casenave}

\maketitle

\begin{abstract}
In an industrial group like Safran, numerical simulations of physical phenomena are integral to most design processes. At Safran's corporate research center, we enhance these processes by developing fast and reliable surrogate models for various physics.
We focus here on two technologies developed in recent years. The first is a physical reduced-order modeling method for non-linear structural mechanics and thermal analysis, used for calculating the lifespan of high-pressure turbine blades and performing heat analysis of high-pressure compressors. The second technology involves learning physics simulations with non-parameterized geometrical variability using classical machine learning tools, such as Gaussian process regression. Finally, we present our contributions to the open-source and open-data community.
\end{abstract}

\section{Introduction}
\label{sec:intro}

\begin{figure*}[h!]
    \centering
    \includegraphics[width=\textwidth]{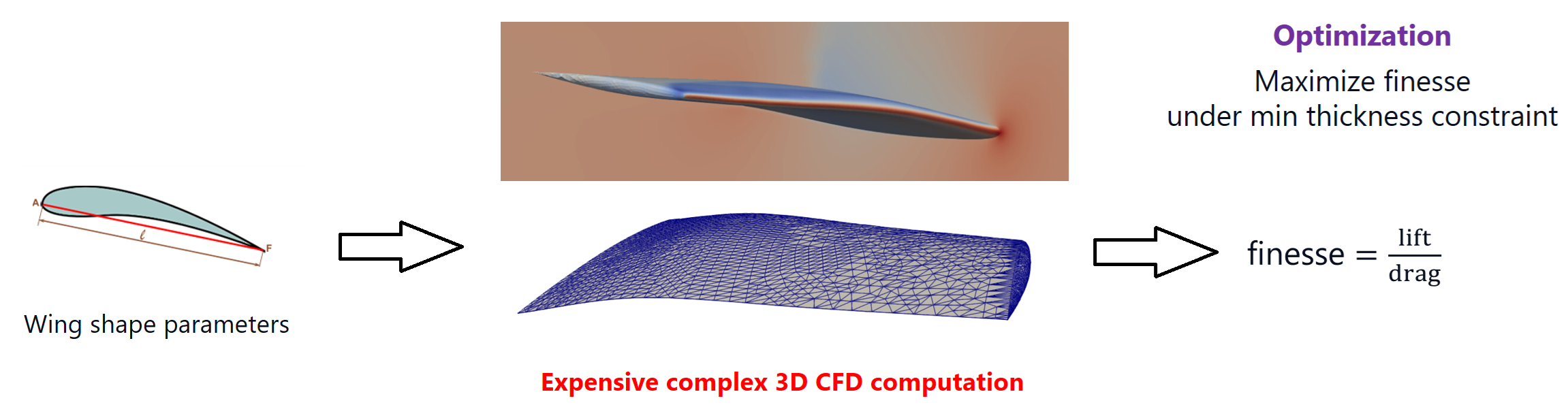}
    \caption{Optimization of the shape of an aircraft wing.}
    \label{fig:wing_optim.png}
\end{figure*}

Consider an operator $F$ that represents a complex numerical problem or an expensive experiment. Industrial design processes often require frequent evaluations of such operators in many-query tasks, such as parametric exploration, optimization, uncertainty quantification, or calibration. These tasks are intractable without modifications or simplifications, making them unsuitable for real-life industrial applications. For instance, optimizing the shape of an aircraft wing to maximize finesse under a minimum thickness constraint involves many iterations of running a complex fluid solver, which is computationally intensive, see Figure~\ref{fig:wing_optim.png}.

A common approach to manage these tasks is by replacing the expensive operator $F$ with a surrogate $f$. This surrogate is created by learning the behavior of $F$ from data. A training database is constructed by selecting inputs using a Design of Experiments (DoE) algorithm and computing corresponding outputs by evaluating $F$. The surrogate, formed as a regressor trained on these input/output pairs, must be accurate enough to replace $F$ in many-query tasks and fast enough to be practical for industrial applications.

For tabular problems where inputs and outputs are low-dimensional vectors ($\mathbb{R}^d \ni x \mapsto y=F(x)\in\mathbb{R}^p$), traditional regression methods like linear or polynomial regression, nearest neighbors, random forests, support vector regression, or Gaussian process regression can be used. This work focuses on learning physics problems, which involve complex, high-dimensional heterogeneous objects like meshes, boundary conditions, fields, or material models, unlike tabular problems. The inputs are the varying components of the physical setting in the many-query task, while the outputs are the predicted components used for decision-making.

For example, the AirfRANS physics learning problem~\cite{bonnet2022airfrans} involves a mesh of the fluid domain, boundary conditions, the Reynolds-Averaged Navier–Stokes equations, and the $k-\omega$ turbulence model. Here, the inputs include the mesh, angle of attack, and inlet velocity, while the outputs are the velocity, pressure, and turbulent velocity fields. The chosen turbulence model being fixed in the training dataset, it is not part of the inputs of the learning problem.

The advantage of physics learning problems over tabular problems with scalar inputs and outputs is the flexibility to change shapes beyond the training parameterization and the ability to modify the post-treatments of predicted solution fields after learning.

This article first discusses physical reduced-order modeling technologies developed at Safran in Section~\ref{sec:PROM}. Then, it introduces a machine learning method for non-parameterized geometrical variability in Section~\ref{sec:ml4phys}. Finally, it presents Safran's initiatives in the open-source and open-data communities in Section~\ref{sec:open}.

The article is a shortened version of the paper submitted to the Proceedings of the 9th European Congress of Mathematics, based on the work presented during the Felix Klein prize lecture on July 16th, 2024.

\section{Physical Reduced-Order Modeling}
\label{sec:PROM}

Physical Reduced-Order Modeling (ROM) is a surrogate modeling technique that uses the underlying physics equations during the exploitation stage.

\subsection{Linear data compression}
\label{sec:POD-ECM}

During the training stage, we generate inputs using DoE algorithms and compute corresponding outputs by evaluating the operator $F$. We then apply a linear dimensionality reduction algorithm to these outputs to construct a reduced-order basis. In the exploitation stage, we use the same algorithms as in the high-dimension problem but restrict the solution to the subspace spanned by the reduced-order basis. For example, if $F$ involves solving a high-dimension initial boundary value problem, we replace the finite element basis with the reduced-order basis to keep using the Galerkin method, expressing the reduced solution as a linear combination of the reduced-order basis elements. Various methods have been proposed to adapt physical ROMs to different problems and simulation methods involved in $F$ evaluations. One such method is the Reduced Basis method~\cite{paterabook}, which starts with a single input/output pair and enriches the data by selecting configurations where the current model makes the largest error, using fast and accurate error bounds~\cite{casenave2014accurate} for accuracy assessments and to trigger enrichment steps.

The speed-up in physical ROMs comes from the reduced-order basis having a much smaller cardinality than the finite element basis, leading to smaller linear systems to solve in the Newton's algorithm during the exploitation stage. However, constructing these small systems can be costly, especially with non-linear problems requiring numerical integration over the mesh. To achieve practical speed-up, precomputations and approximations (called hyperreduction) are performed during the training stage to reduce computational effort in the exploitation stage, see~\cite[Section 3.2]{casenave2020nonintrusivemech} for more details, and~\cite[Section 1]{casenave2020nonintrusivetherm} for references on hyperreduction methods.

We focus on non-linear structural mechanics and non-linear transient heat problems, utilizing snapshot Proper Orthogonal Decomposition (POD) for dimensionality reduction, which is parallelizable with domain decomposition. Additionally, we use a localized Empirical Cubature Method (ECM), which derives reduced quadrature formulae independently for each subdomain. This involves solving an optimization problem to minimize the number of non-zero elements of the quadrature weights vector, using a greedy approach called non-negative orthogonal matching pursuit. This ensures the reduced operators retain the spectral properties of $F$, maintaining symmetry and positive definiteness of the tangent operator in the Newton's algorithm.

For non-linear structural mechanics and transient heat problems, we have demonstrated significant speed-ups, see respectively~\cite{casenave2020nonintrusivemech} and~\cite{casenave2020nonintrusivetherm}. For instance, our physical ROM reduced the computation time for lifetime prediction of 3D elastoviscoplastic high-pressure turbine blades from weeks to hours, and accelerated transient heat analysis of high-pressure compressors with complex boundary conditions, achieving errors compatible with industrial requirements.

In structural mechanics applications, quantities of interest often include dual quantities, such as plastic strain or Von Mises stress. To reconstruct dual quantities over the entire mesh, the Gappy-POD method is used. By enforcing the reduced quadrature scheme to include well-chosen quadrature points, the Gappy-POD reconstruction can be made well-posed, see~\cite[Proposition~1]{casenave2020nonintrusivemech}. Our numerical experiments revealed a strong correlation between the prediction error on dual quantities and the ROM-Gappy-POD residual, which is defined as the error, on the reduced quadrature points, between the online computation of the corresponding quantity by the constitutive law solver and the Gappy-POD reconstruction. Thus, we introduced~\cite{casenave2019error} an error indicator, which is a regressor trained to predict the error from the ROM-Gappy-POD residual. This indicator successfully maintains the error below a chosen threshold by triggering enrichment steps.

Linear dimensionality reduction enables the application of the Galerkin method during the exploitation stage. However, this restriction can sometimes result in a reduced-order basis with a cardinality too large to achieve significant speed-up, even with hyperreduction. These situations, referred to as poorly reducible problems, will be addressed in the following section.

\subsection{Piecewise linear data compression}

The reducibility of the approximation of $F$ by $f$ linear dimensionality reduction is evaluated by the rate of decay of the Kolmogorov $N$-width of the output solutions with respect to the dimension $N$ of an optimal approximation vector space. The Kolmogorov $N$-width measures the worst approximation error of a set of outputs by the best $N$-dimensional subspace. If the decay rate is fast, linear data compression methods are effective. If not, the reduced-order basis becomes too large, making the reduced problems slow to solve, indicating a poorly reducible problem. This challenge is addressed in Thomas Daniel's Ph.D. work~\cite{danie2021}.

\begin{figure}[h!]
    \centering
    \includegraphics[width=0.38\textwidth]{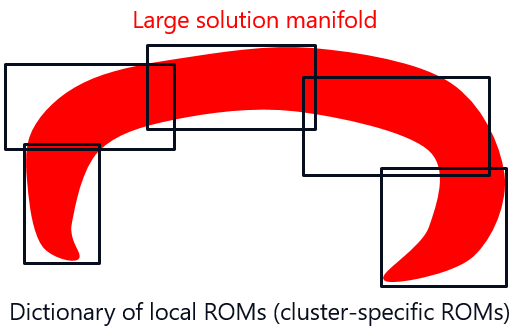}
    \caption{Piecewise linear data compression}
    \label{fig:piecewise_lin_data_comp}
\end{figure}

\begin{figure*}[h!]
    \centering
    \includegraphics[width=\textwidth]{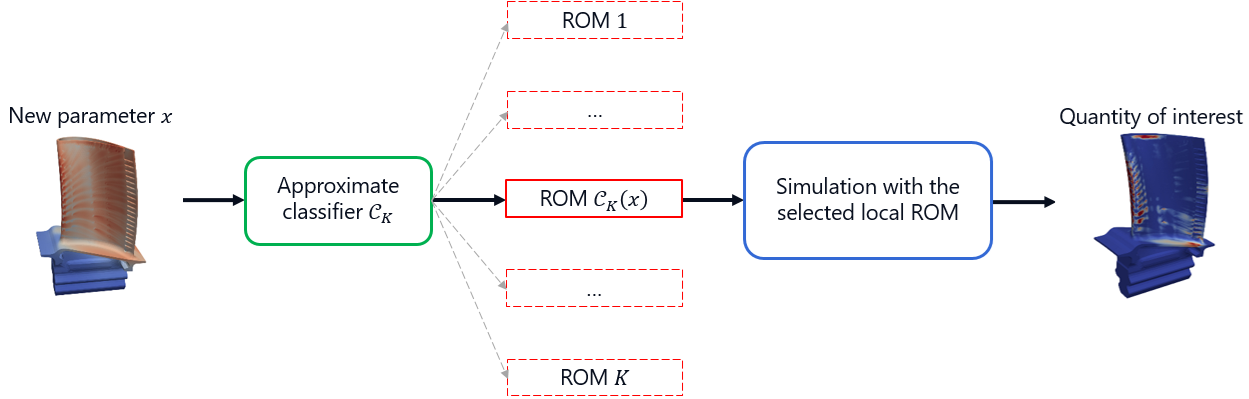}
    \caption{Local ROM recommendation.}
    \label{fig:classifier_ROM}
\end{figure*}

To tackle poorly reducible problems, we focus on creating a dictionary of local physical ROMs using piecewise linear data compression, by clustering the outputs fields, see Figure~\ref{fig:piecewise_lin_data_comp}.
Consider the problem of predicting the lifetime of 3D elastoviscoplastic high-pressure turbine blades, now with the temperature loading field as input. To use a dictionary of local ROMs, we need to know in advance in which cluster of the outputs the solution is located, for us to use the adequate local ROM. To achieve this, we train a classifier to select the correct local ROM based on the 3D input temperature field, similar to how classifiers in computer graphics recognize objects in images, see Figure~\ref{fig:classifier_ROM}. Such workflows are proposed in~\cite{daniel2020model}. Unlike image classifiers, our situation involves few high-dimensional inputs. To improve the classifier training, we developed specialized feature selection and data augmentation techniques. This approach ensures accurate local ROM selection and effective handling of complex, high-dimensional input data, see~\cite{daniel2021data}.

Partitioning outputs using the $L^2$ norm in the ambient solution space has been explored in the literature, but this does not guarantee low-dimensional approximation subspaces. In~\cite[Property 4.13]{daniel2022physics}, we demonstrate that the partitions of outputs minimizing the k-medoid cost function with sine dissimilarity are are exactly the minimizers of the sum of a local variant of the Kolmogorov 1-width, weighted by the partition element's volume. This method, which considers the relative angles between outputs, offers a practical algorithm for optimally partitioning outputs for our purposes.

\begin{figure*}[h!]
    \centering
    \includegraphics[width=0.7\textwidth]{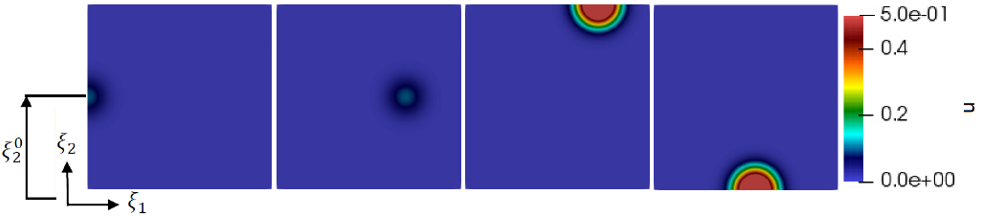}
    \caption{Some output fields for the 2D advection problem~\cite{daniel2022physics}.}
    \label{fig:transport_2D_visu}
\end{figure*}

Consider a 2D advection problem where Gaussians with small (0.1) and large (1) amplitudes and various vertical positions $\xi_2^0$ move from left to right, see Figure~\ref{fig:transport_2D_visu}. 
We illustrate the clustering of these outputs using multidimensional scaling (MDS) with five clusters, comparing $L^2$ norm and sine dissimilarity. With the $L^2$ norm, all outputs with the small amplitude (0.1) are tightly packed at the center, likely grouping into the same cluster. This cluster contains all the independent directions of the outputs, hence has a Kolmogorov $N$-width not smaller than the one of the set of all outputs. Conversely, with sine dissimilarity, three connected components emerge, corresponding to the three considered $\xi_2^0$ values. At each position are located two outputs varying only by magnitude, indicating zero sine dissimilarity for outputs differing only in magnitude.

\begin{figure*}[h!]
    \centering
    \includegraphics[width=0.6\textwidth]{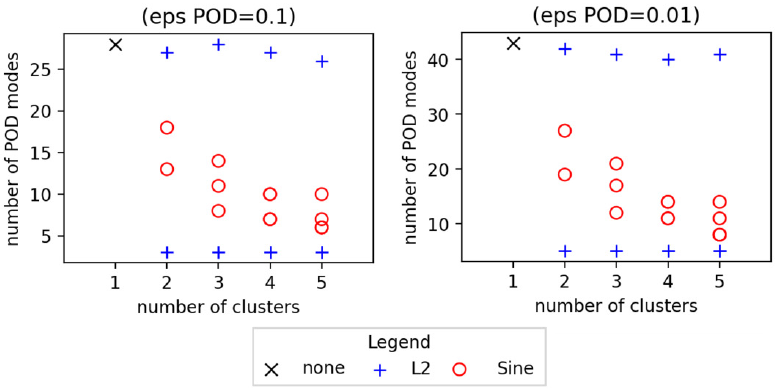}
    \caption{Number of POD modes for each local basis with respect to the number of clusters, for various accuracy criterions of the POD truncature, applied to the 2D advection problem~\cite{daniel2022physics}.}
    \label{fig:transport_2D_comp}
\end{figure*}

In Figure~\ref{fig:transport_2D_comp}, we compare the effectiveness of clustering outputs using the $L^2$ norm versus sine dissimilarity in constructing local ROMs of small dimensions. When using the $L^2$ norm, one cluster ends up with the same dimension as a global ROM, specifically the cluster containing all small magnitude outputs. This indicates that no matter how many clusters are chosen, one cluster remains as poorly reducible as the entire set of outputs. Conversely, when using sine dissimilarity, the dimensions of all local ROMs decrease as the number of clusters increases. This demonstrates that clustering outputs using sine dissimilarity is beneficial, as it successfully reduces the ROM dimensions. Therefore, the choice of dissimilarity in output clustering is crucial, and the sine dissimilarity proves to be effective in this context.

In~\cite{daniel2022uncertainty}, a dictionary of local ROMs is developed to quantify the uncertainty of the plastic strain and Von Mises stress in high-pressure turbine blades with respect to input temperature loading. First, we parameterize the variability of the input temperature field and constructe two DoEs: one MaxProj of size 80 and one Sobol' of size 120. We evaluate the operator $F$ for these 200 configurations, taking 7 days and 9 hours on 48 computer cores. Using the sine dissimilarity with the PAM k-medoid clustering algorithm, we partition the MaxProj DoE outputs into two clusters in 5 minutes. Next, we construct the two corresponding local ROMs using the POD-ECM technique presented in Section~\ref{sec:POD-ECM}, which takes 5 hours on 24 computer cores. We label the 120 Sobol' DoE outputs by computing the sine dissimilarity with the two identified medoids in 5 minutes. Our feature selection algorithm takes 16 minutes on 280 computer cores, followed by training a logistic regression classifier with elastic-net regularization in 1 minute. This process results in a ROM dictionary of two local ROMs and a classifier capable of recommending the appropriate local ROM based on input temperature loading, making our surrogate $f$ ready for use. Finally, we conduct probability density estimation of the plastic strain and Von Mises stress in areas of interest on the turbine blade using 1,008 Monte Carlo evaluations of $f$. This yields relative errors of 1-2\% and a speed-up factor of 600.

\subsection{Limitations}

Despite the successful application of physical ROMs to many equations and configurations, two significant challenges persist. First, non-parameterized geometrical variations make it difficult to compare outputs supported on different meshes. Second, global non-linear data compression prevents the use of the Galerkin method during the exploitation stage. Some authors have developed solutions to address these limitations while still constructing a surrogate $f$ that involves solving the physics equations. However, these approaches often involve complex algorithms that can hinder practical speed-ups. 

In the following section, we propose a surrogate $f$ that incorporates non-linear data compression and manages non-parameterized geometrical variations, but does not rely on solving the physics equations during the exploitation stage.

\section{Learning physical problems with non-parameterized geometrical variability}
\label{sec:ml4phys}

\subsection{Mesh Morphing Gaussian Process}

In our recent work~\cite{casenave2024mmgp}, we introduce the Mesh Morphing Gaussian Process (MMGP) technique, which integrates four key components: (i) mesh morphing, (ii) finite element interpolation, (iii) dimensionality reduction, and (iv) Gaussian process regression. This approach is designed to learn solutions to partial differential equations (PDEs) that involve geometric variations not explicitly parameterized.

\begin{figure*}[h!]
\centering
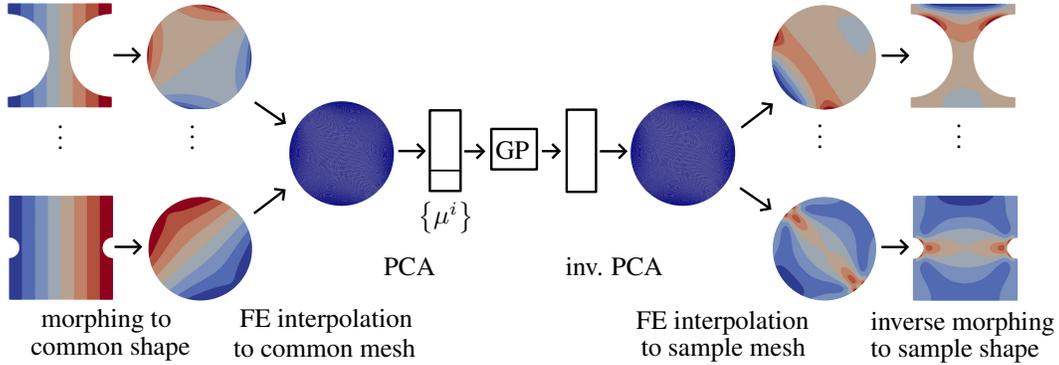
\caption{Illustration of the MMGP inference workflow for the prediction of an output field of interest~\cite{casenave2024mmgp}.}
\label{fig:mmgp}
\end{figure*}

Figure~\ref{fig:mmgp} illustrates the MMGP inference process for predicting a field of interest. The figure should be read from left to right: the left side shows various geometries corresponding to different samples, while the right side displays the target fields to be predicted on these geometries.
For non-parameterized meshes, we perform a shape embedding process to convert meshes into learnable objects, by considering the coordinates of the mesh vertices as continuous fields. The left column of Figure~\ref{fig:mmgp} depicts the continuous field of the x-component of the coordinates, showing vertical iso-values. We then apply a deterministic morphing process to transform input meshes into a common shape, such as the unit disk shown in Figure~\ref{fig:mmgp}. Each sample is converted into a mesh of the unit disk, modifying the coordinate fields according to the input mesh shapes. We then select a common mesh of the unit disk and project the coordinate fields from each morphed mesh onto this common mesh using finite element interpolation. This process results in all coordinate fields being represented on the same mesh, allowing us to use classical dimensionality reduction techniques, like Principal Component Analysis (PCA), to generate low-dimensional vectors for the meshes. If scalar inputs are part of the learning problem, these are appended at the end of these vectors. 

Similarly, a field embedding process is needed to convert variable-sized output fields into learnable objects. We apply the same transformations used for coordinate fields: morphing onto a common shape, finite interpolation onto a common mesh, and PCA compression. This process results in low-dimensional vectors representing our output fields of interest.

These deterministic pre-treatments effectively reduce the physics learning problem to a low-dimensional input/output tabular regression problem. The morphing and related steps act as the dimensionality reduction phase, which, while possibly highly non-linear, simplifies the machine learning task by converting large, variable-dimensional objects into smaller manageable ones. The deterministic nature of these pre-treatments is essential as it simplifies the machine learning stage and avoids the complexity of dealing with large and variable-dimensional data directly. 

We use Gaussian process regression for its high accuracy and ability to estimate predictive uncertainties effectively. For further details and numerical experiments demonstrating the efficacy of MMGP, including comparisons with recent deep learning technologies like MeshGraphNet and U-net graph convolutional networks, we refer to~\cite{casenave2024mmgp}.

\subsection{ML4PhySim competition}
\label{sec:ml4physim}

We participated in a competition~\cite{irtchallenge} organized by IRT SystemX and partners, where we employed our MMGP technology. The challenge involved predicting pressure, velocity, and turbulent viscosity fields on the AirfRANS dataset~\cite{bonnet2022airfrans}, introduced in Section~\ref{sec:intro}.

\begin{figure*}[h!]
    \centering
    \includegraphics[width=0.7\textwidth]{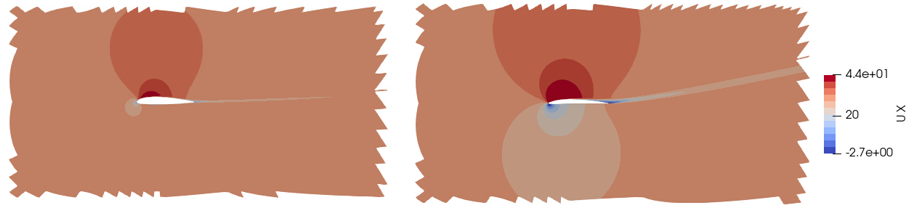}
    \caption{Illustration of the horizontal velocity for two samples of the AirfRANS dataset: the external boundaries are very irregular.}
    \label{fig:challenge_samples}
\end{figure*}

To achieve effective morphing, we first extended all 2D meshes to a common boundary box to handle irregular external boundaries, as shown in Figure~\ref{fig:challenge_samples}.

We select the geometry of the first pretreated mesh from the training set as the common shape. We specify the morphing of chosen control points—on the intrado and extrado of the airfoil, the wake line, and the external boundary-and use a Radial Basis Function (RBF) morphing algorithm to compute the transformation of all other points. We project all fields, including input coordinates and output fields, onto the common mesh (the pretreated and morphed mesh from the first training sample). For dimensionality reduction, we used snapshot Proper Orthogonal Decomposition (POD) as detailed in~\cite[Annex C]{casenave2024mmgp}. The snapshot POD approach ensures that objects reconstructed from reduced representations preserve zero linear relationships, which is critical for maintaining zero boundary conditions, such as those for the velocity field at the airfoil surface. This helps to ensure that the zero boundary condition is respected in the predicted fields. We train separate Gaussian Processes for each generalized coordinate associated with the snapshot-POD of each output field. The Gaussian Process kernel combines a constant term with an RBF kernel and a white noise term.

\begin{figure*}[h!]
    \centering
    \includegraphics[width=\textwidth]{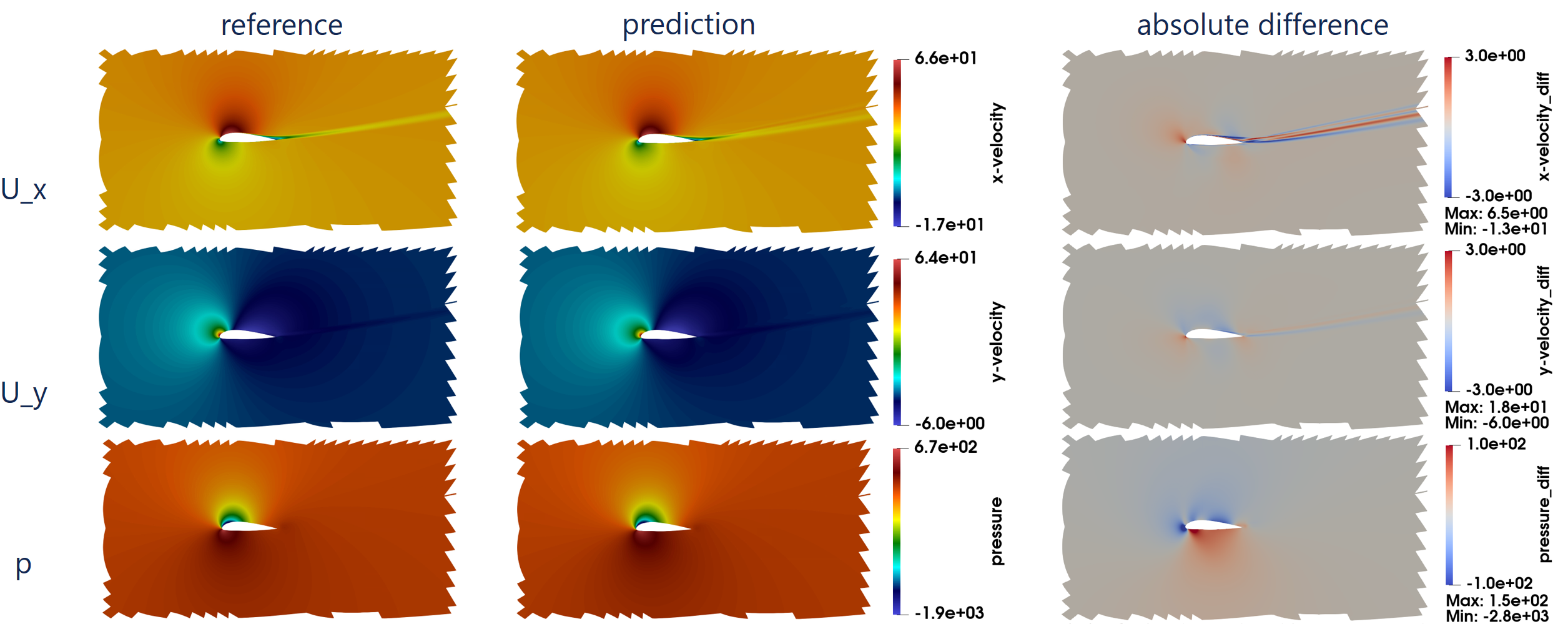}
    \caption{Illustration of the result on a test sample.}
    \label{fig:1_test_ood_sample}
\end{figure*}

Some results are illustrated in Figure~\ref{fig:1_test_ood_sample}.

\begin{figure*}[h!]
    \centering
    \includegraphics[width=\textwidth]{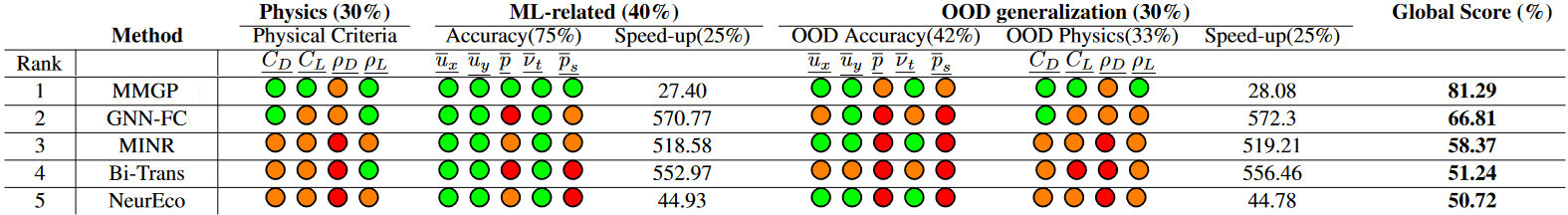}
    \caption{ML4PhySim competition final leaderboard~\cite{irtchallenge2}.}
    \label{fig:challenge_leaderboard}
\end{figure*}

The competition score, ranging from 0 to 100, evaluates both the speed and accuracy of the surrogate model, including scalar quantities derived from the output fields (such as lift and drag coefficients) and their Spearman's rank correlation coefficients, tested on two test sets, one of them being out-of-distribution. Figure~\ref{fig:challenge_leaderboard} shows that our method won first place, outperforming competitors who used advanced deep learning techniques, despite our approach using simpler and more classical tools.

\subsection{Outlooks}
\label{sec:mmgp_conclusion}

MMGP is an effective method for approximating field and scalar outputs in complex physics problems with non-parameterized geometrical variability. It combines mesh morphing with finite element interpolation and uses shape embedding by reducing the dimensions of coordinate fields, simplifying the machine learning task. This approach enables efficient Gaussian process regression on reduced-dimensional data.

The method can handle large meshes, can be trained on CPU, is interpretable, and provides accurate predictions and predictive uncertainties. It is especially efficient for industrial part design, where configurations often have low intrinsic dimensions.

However, MMGP has some limitations: the morphing process must be customized for each use case, it assumes a fixed mesh topology, and speed-up is limited by the need for morphing and interpolation during inference. Future work will focus on addressing these issues by developing a generic morphing method, creating efficient approximate morphing procedures for inference~\cite{kabalan2024elasticity}, and deriving an optimal morphing strategy to minimize dimensionality further.

\section{Open-source and open-data}
\label{sec:open}

Since the inauguration of the Safran corporate research center in 2015, we have implemented processes to make our research codes open-source and non-confidential datasets available as open-data. Our commitment to accessibility, quality, and robustness involves using DevOps tools and practices, such as git for version control, code reviews, automated testing, deployment pipelines, and extensive documentation. We opt for permissive licenses like MIT or BSD3 for our code and CC-BY-SA for datasets. Our codes are hosted on GitLab, documentation on Read the Docs, and packages are distributed via conda-forge, when possible.

\subsection{Open-source codes}

One recent development is the PLAID (Physics Learning AI Datamodel) library, which formalizes physics learning problems and manages complex datasets. It leverages the CGNS standard~\cite{cgns} for data handling, promoting consistency and ease of data sharing among practitioners.

The codes developed for the methods presented in the previous sections are open-source. Notably, Muscat~\cite{bordeu2023basictools, muscat} is a library for mesh processing in finite element computations. PLAID \cite{plaid} is used as a foundation for our physics learning developments. Physical ROM technologies are supported by two libraries: mordicus~\cite{mordicus}, co-developed with the MOR\_DICUS consortium, and genericROM~\cite{genericROM}, developed by Safran. The implementation of genericROM is non-intrusive, parallel with distributed memory, and can handle non-parameterized variability, see~\cite{ryckelynck2024resources} for details on these notions and a description of mordicus and genericROM.
Additionally, MMGP, discussed in Section~\ref{sec:ml4phys}, is also available as open-source~\cite{mmgp}. For tabular learning problems, Safran's open-source platform Lagun~\cite{lagun} supports a range of tasks including design of experiments, surrogate modeling, sensitivity analysis, and optimization. 

\subsection{Open-data datasets}

We are in the process of making our non-confidential datasets available as open-data in the PLAID format, facilitating their use by the community. These datasets are shared through Zenodo and Hugging Face. Specifically, the three datasets used in our recent comparison of MMGP with advanced deep learning methods in~\cite{casenave2024mmgp} are accessible on Zenodo~\cite{casenave_2023_zenodo, roynard_2024_zenodo, roynard_2023_zenodo} and Hugging Face~\cite{casenave_2023_huggingface, roynard_2024_huggingface, roynard_2023_huggingface}.

\begin{figure*}[h!]
    \centering
    \includegraphics[width=\textwidth]{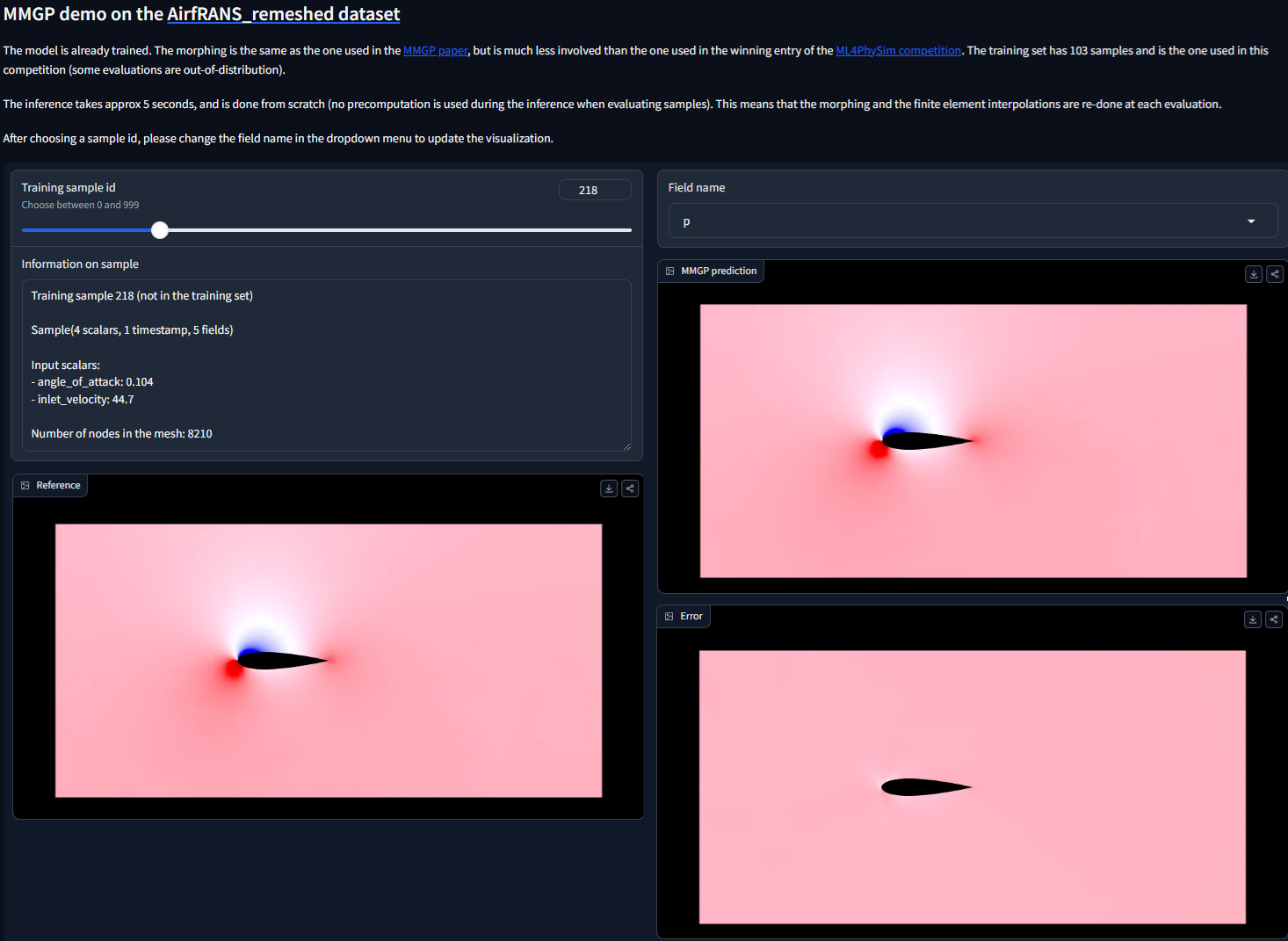}
    \caption{Hugginface space hosting a demo of our method MMGP applied on the physics learning problem of the ML4PhySim competition.}
    \label{fig:mmgp_demo}
\end{figure*}

Additionally, we have distributed a demo of our MMGP method on Hugging Face's Spaces platform, which was applied to the ML4PhySim competition problem presented in Section~\ref{sec:ml4physim}, see Figure~\ref{fig:mmgp_demo}.
In this demo, users can select a sample ID, which triggers the MMGP surrogate model, which involves: morphing, finite element interpolation, projection onto snapshot-POD modes, inference using a Gaussian process regressor, and finally, inverse snapshot-POD and finite element interpolation. Predicted fields, reference data, and errors are plotted. Note that for real-time performance, the demo uses a remeshed dataset with meshes reduced from approximately 180,000 nodes to around 8,000 nodes.

\bibliography{main}

\begin{thebibliography}{10}

\bibitem{cgns}
Cfd general notation system (cgns).
\newblock \url{https://cgns.github.io}.
\newblock Accessed: 2024-08-06.

\bibitem{bonnet2022airfrans}
F.~Bonnet, J.~Mazari, P.~Cinnella, and P.~Gallinari.
\newblock Airfrans: High fidelity computational fluid dynamics dataset for approximating reynolds-averaged navier--stokes solutions.
\newblock {\em Advances in Neural Information Processing Systems}, 35:23463--23478, 2022.

\bibitem{bordeu2023basictools}
F.~Bordeu, F.~Casenave, and J.~Cortial.
\newblock Basictools: a numerical simulation toolbox.
\newblock {\em Journal of Open Source Software}, 8(86):5142, 2023.

\bibitem{casenave2019error}
F.~Casenave and N.~Akkari.
\newblock An error indicator-based adaptive reduced order model for nonlinear structural mechanics—application to high-pressure turbine blades.
\newblock {\em Mathematical and Computational Applications}, 24(2), 2019.

\bibitem{casenave2020nonintrusivemech}
F.~Casenave, N.~Akkari, F.~Bordeu, C.~Rey, and D.~Ryckelynck.
\newblock A nonintrusive distributed reduced-order modeling framework for nonlinear structural mechanics—application to elastoviscoplastic computations.
\newblock {\em International Journal for Numerical Methods in Engineering}, 121(1):32--53, 2020.

\bibitem{casenave2014accurate}
F.~Casenave, A.~Ern, and T.~Leli{\`e}vre.
\newblock Accurate and online-efficient evaluation of the a posteriori error bound in the reduced basis method.
\newblock {\em ESAIM: Mathematical Modelling and Numerical Analysis}, 48(1):207--229, 2014.

\bibitem{casenave2020nonintrusivetherm}
F.~Casenave, A.~Gariah, C.~Rey, and F.~Feyel.
\newblock A nonintrusive reduced order model for nonlinear transient thermal problems with nonparametrized variability.
\newblock {\em Advanced Modeling and Simulation in Engineering Sciences}, 7:1--19, 2020.

\bibitem{casenave_2023_zenodo}
F.~Casenave, X.~Roynard, and B.~Staber.
\newblock {Tensile2d: 2D quasistatic non-linear structural mechanics solutions, under geometrical variations}, November 2023.

\bibitem{casenave_2023_huggingface}
F.~Casenave, X.~Roynard, and B.~Staber.
\newblock {Tensile2d: 2D quasistatic non-linear structural mechanics solutions, under geometrical variations}, 2024.

\bibitem{casenave2024mmgp}
F.~Casenave, B.~Staber, and X.~Roynard.
\newblock {MMGP}: a {M}esh {M}orphing {G}aussian {P}rocess-based machine learning method for regression of physical problems under nonparametrized geometrical variability.
\newblock {\em Advances in Neural Information Processing Systems}, 36, 2024.

\bibitem{mordicus}
MOR\_DICUS consortium.
\newblock Mordicus.
\newblock code: \url{https://gitlab.com/mor_dicus/mordicus}, doc: \url{https://mordicus.readthedocs.io}.
\newblock Accessed: 2024-08-07.

\bibitem{danie2021}
T.~Daniel.
\newblock {\em Machine learning for nonlinear model order reduction}.
\newblock PhD thesis, 2021.
\newblock Thèse de doctorat dirigée par Ryckelynck, David Mécanique Université Paris sciences et lettres 2021.

\bibitem{daniel2022physics}
T.~Daniel, F.~Casenave, N.~Akkari, A.~Ketata, and D.~Ryckelynck.
\newblock Physics-informed cluster analysis and a priori efficiency criterion for the construction of local reduced-order bases.
\newblock {\em Journal of Computational Physics}, 458:111120, 2022.

\bibitem{daniel2020model}
T.~Daniel, F.~Casenave, N.~Akkari, and D.~Ryckelynck.
\newblock Model order reduction assisted by deep neural networks (rom-net).
\newblock {\em Advanced Modeling and Simulation in Engineering Sciences}, 7:1--27, 2020.

\bibitem{daniel2021data}
T.~Daniel, F.~Casenave, N.~Akkari, and D.~Ryckelynck.
\newblock Data augmentation and feature selection for automatic model recommendation in computational physics.
\newblock {\em Mathematical and Computational Applications}, 26(1):17, 2021.

\bibitem{daniel2022uncertainty}
T.~Daniel, F.~Casenave, N.~Akkari, D.~Ryckelynck, and C.~Rey.
\newblock Uncertainty quantification for industrial numerical simulation using dictionaries of reduced order models.
\newblock {\em Mechanics \& Industry}, 23:3, 2022.

\bibitem{irtchallenge}
IRT~SystemX et~al.
\newblock Ml4physim challenge.
\newblock \url{https://www.codabench.org/competitions/1534}.
\newblock Accessed: 2024-08-07.

\bibitem{irtchallenge2}
IRT~SystemX et~al.
\newblock Neurips 2024 - ml4cfd competition.
\newblock \url{https://www.codabench.org/competitions/3282}.
\newblock Accessed: 2024-08-07.

\bibitem{kabalan2024elasticity}
A.~Kabalan, F.~Casenave, F.~Bordeu, V.~Ehrlacher, and A.~Ern.
\newblock An elasticity-based mesh morphing technique with application to reduced-order modeling.
\newblock {\em arXiv preprint arXiv:2407.02433}, 2024.

\bibitem{paterabook}
A.~T. Patera and G.~Rozza.
\newblock {\em Reduced Basis Approximation and A Posteriori Error Estimation for Parametrized Partial Differential Equations}.
\newblock MIT Pappalardo Graduate Monographs in Mechanical Engineering, 2007.

\bibitem{roynard_2023_zenodo}
X.~Roynard, F.~Casenave, and B.~Staber.
\newblock {Rotor37: a 3D CFD RANS dataset, under geometrical variations of a compressor blade}, November 2023.

\bibitem{roynard_2024_zenodo}
X.~Roynard, F.~Casenave, and B.~Staber.
\newblock Airfrans\_original, June 2024.

\bibitem{roynard_2024_huggingface}
X.~Roynard, F.~Casenave, and B.~Staber.
\newblock Airfrans\_original, 2024.

\bibitem{roynard_2023_huggingface}
X.~Roynard, F.~Casenave, and B.~Staber.
\newblock {Rotor37: a 3D CFD RANS dataset, under geometrical variations of a compressor blade}, 2024.

\bibitem{ryckelynck2024resources}
D.~Ryckelynck, F.~Casenave, and N.~Akkari.
\newblock Resources: Software and tutorials.
\newblock In {\em Manifold Learning: Model Reduction in Engineering}, pages 53--69. Springer, 2024.

\bibitem{genericROM}
Safran.
\newblock {genericROM}.
\newblock code: \url{https://gitlab.com/drti/genericrom}, doc: \url{https://genericrom.readthedocs.io}.
\newblock Accessed: 2024-08-07.

\bibitem{lagun}
Safran.
\newblock Lagun.
\newblock code: \url{https://gitlab.com/drti/lagun}, doc: \url{https://lagun.readthedocs.io}.
\newblock Accessed: 2024-08-07.

\bibitem{mmgp}
Safran.
\newblock Mmgp.
\newblock code: \url{https://gitlab.com/drti/mmgp}, doc: \url{https://mmgp.readthedocs.io}.
\newblock Accessed: 2024-08-07.

\bibitem{muscat}
Safran.
\newblock Muscat.
\newblock code: \url{https://gitlab.com/drti/muscat}, doc: \url{https://muscat.readthedocs.io}.
\newblock Accessed: 2024-08-07.

\bibitem{plaid}
Safran.
\newblock Plaid.
\newblock code: \url{https://gitlab.com/drti/plaid}, doc: \url{https://plaid-lib.readthedocs.io}.
\newblock Accessed: 2024-08-07.

\end{thebibliography}
\bibliographystyle{plain}

\end{document}